\documentclass[12pt]{amsart}
\usepackage{amssymb,amscd,amsthm,verbatim,amsmath}

\makeatletter \headheight=8pt     \topmargin=10pt
\newtheorem{proposition}{Proposition}[section]
\newtheorem{theorem}[proposition]{Theorem}
\newtheorem{corollary}[proposition]{Corollary}
\newtheorem{lemma}[proposition]{Lemma}
\newtheorem{conjecture}[proposition]{Conjecture}

\newtheorem{remark}[proposition]{Remark}

\newtheorem{definition}[proposition]{Definition}

\newcommand\supp{\operatorname{sp}}

\newcommand{\C}{\mathcal{C}}
\newcommand{\D}{\mathcal{D}}

\newcommand{\F}{\mathcal{F}}
\renewcommand{\L}{\mathcal{L}}
\newcommand{\T}{\mathcal{T}}
\newcommand{\rank}{{\mathrm{rank}}}
\newcommand{\link}{{\mathrm{link}}}

\begin{document}

\title[Uniform affine oriented matroids]
{The bounded complex of a uniform affine oriented matroid is a
ball}

\author{Xun Dong}
\email{xundong@math.miami.edu}

\address{Department of Mathematics, University of Miami, Coral Gables, FL 33124}


\begin{abstract}
Zaslavsky~\cite{Za} conjectures that the bounded complex of a
simple hyperplane arrangement is homeomorphic to a ball. We prove
this conjecture for the more general uniform affine oriented
matroids.
\end{abstract}

\maketitle

\section{Introduction}

A hyperplane arrangement in a Euclidean space partitions the space
into faces. The bounded faces are those bounded in the usual
metric sense. The collection of all bounded faces is a polyhedral
complex called the {\em bounded complex} of the hyperplane
arrangement. More generally, affine oriented matroids have a
topological model as arrangements of {\em pseudo-hyperplanes},
each obtained from a flat hyperplane by tame topological
deformation. Since we can still talk about faces and metric in a
pseudo-hyperplane arrangement, the {\em bounded complex} of an
affine oriented matroid can be defined as the regular cell complex
consisting of all bounded faces. Seemingly a simple object, the
topology of the bounded complex is still not completely
understood.

Study on the bounded complex goes back to Zaslavsky's 1975
paper~\cite{Za}, in which he proves several face counting formulas
for hyperplane arrangements. One of the formula is for the number
of bounded regions of a hyperplane arrangement. The formula can be
proven using M\"obius inversion provided that one knows that the
Euler characteristic of the bounded complex is one. This prompts
Zaslavsky~\cite{Za} to conjecture that the bounded complex is
contractible. This conjecture was proven in Ziegler~\cite{Zi} for
hyperplane arrangements, and in Bj\"orner and Ziegler~\cite{BZ}
for affine oriented matroids. In fact, Zaslavsky~\cite{Za} also
conjectures that the bounded complex is star-convex, that is,
there is a ``center point'' from which all other points in the
bounded complex can be seen. However this is false (see, eg.,
\cite[Exercise 4.29]{BLSWZ}). Here I should remark that some of
Zaslavsky's formulas were independently discovered by Las
Vergnas~\cite{L}. The reader may consult~\cite[\S 4.6]{BLSWZ} for
an account of relevant formulas and their histories.

Zaslavsky~\cite{Za} also proves that the bounded complex of a
hyperplane arrangement is pure, that is, all the maximal bounded
faces have the same dimension. However his proof does not
generalize to affine oriented matroids. In a recent paper~\cite{D}
it is shown that the bounded complex of an affine oriented matroid
is pure by using {\em covector axioms}. It is also shown
in~\cite{D} that the bounded complex is {\em collapsible}. A
collapsible complex is contractible, but not vice versa. The
collapsibility of the bounded complex will play a crucial role in
this paper.

In general the bounded complex of a hyperplane arrangement is not
necessarily a (closed) ball. For instance, let us consider the
four lines defined by equations $x=0, y=0, x+y=1$ and $x+y=-1$
respectively in a plane. The bounded complex of this line
arrangement consists of two triangles joined at a vertex. However
if the hyperplanes are in general position (such an arrangement is
called {\em simple}), then it is intuitively plausible that the
bounded complex should be a ball. This was conjectured to be the
case in Zaslavsky~\cite{Za} (see also Stanley~\cite{St}). The main
objective of this paper is to prove this conjecture. The rough
idea of the proof is as follows. It is known that a collapsible
{\em piecewise-linear ($PL$) manifold} is a $PL$ ball (however a
contractible $PL$ manifold is not necessarily a ball, see
Mazur~\cite{Ma} for one of the first such examples). Therefore it
is sufficient to show that the bounded complex is a $PL$ manifold
since it is known to be collapsible. This will be accomplished by
showing that the link of a vertex in the order complex
$\Delta(\L^{++})$ is either a $PL$ ball or a $PL$ sphere. The
proof works for the more general uniform affine oriented matroids.

The paper is organized as follows. Notations and preliminary facts
about affine oriented matroids are reviewed in Section 2. Section
3 outlines the proof of Zaslavskey's conjecture. In Section 4, the
link of a vertex in the order complex $\Delta(\L^{++})$ gets a
detailed study, and the proof of the conjecture is completed.
Finally Section 5 summarizes some further open questions.

\section{Affine oriented matroids}

Let us start with a quick review of the necessary definitions and
terminology for oriented matroids. We mostly follow Section 4.1
of~\cite{BLSWZ}. The only difference is the notation for the
support of a sign vector.

Let $E$ be a finite set and consider the sign vectors $X, Y \in
\{+,-,0\}^E$. The {\em support} of a vector $X$ is $\supp(X)=\{e
\in E\ :\ X_e \neq 0\}$; its {\em zero set} is
$$
z(X)=E\setminus \supp(X)=\{e\in E\ :\ X_e = 0\}.
$$
The {\em opposite} of a vector $X$ is $-X$, defined by
$$
(-X)_e = \left\{%
\begin{array}{ll}
    -, & \hbox{if $X_e=+$;} \\
    +, & \hbox{if $X_e=-$;} \\
    0, & \hbox{if $X_e=0$.} \\
\end{array}%
\right.
$$
The {\em zero vector} is 0, with $0_e=0$ for all $e\in E$. The
{\em composition} of two vectors $X$ and $Y$ is $X\circ Y$,
defined by
$$
(X\circ Y)_e=\left\{%
\begin{array}{ll}
    X_e, & \hbox{if $X_e\neq 0$;} \\
    Y_e, & \hbox{otherwise.} \\
\end{array}%
\right.
$$
The {\em separation set} of $X$ and $Y$ is $S(X,Y)=\{e\in E\ :\
X_e=-Y_e\neq0\}$. Notice that $X\circ Y=Y\circ X$ if and only if
$S(X,Y)=\emptyset$, in which case we say that $X$ and $Y$ are {\em
conformal}. We are now ready for the definition of oriented
matroids in terms of covectors (see~\cite[4.1.1]{BLSWZ}).

\begin{definition}[Covector Axioms]\rm
An oriented matroid is a pair $(E,\L)$, where $E$ is a finite set
and $\L\subseteq \{+,-,0\}^E$ is the set of {\em covectors}
satisfying:
\begin{enumerate}
    \item[(L0)] $0\in \L$;
    \item[(L1)] $X\in \L$ implies that $-X\in\L$;
    \item[(L2)] $X,Y\in\L$ implies that $X\circ Y \in \L$;
    \item[(L3)] if $X,Y \in \L$ and $e\in S(X,Y)$ then there
    exists $Z\in \L$ such that $Z_e=0$ and $Z_f=(X\circ
    Y)_f=(Y\circ X)_f$ for all $f\notin S(X,Y)$.
\end{enumerate}
\end{definition}

Let $\leq$ be the partial order on the set $\{+,-,0\}$ defined by
$0<+$ and $0<-$, with $+$ and $-$ incomparable. This induces a
product partial order on $\{+,-,0\}^E$. Thus $Y\leq X$ if and only
if $Y_e\in\{0,X_e\}$ for all $e\in E$. As a subset of
$\{+,-,0\}^E$ the set of covectors $\L$ has an induced partial
order with bottom element 0. Let $\widehat{\L}$ denote the poset
$\L$ with a top element $\hat{1}$ adjoined. Then $\widehat{\L}$ is
a lattice called the {\em big face lattice} of $(E,\L)$. The join
in $\widehat{\L}$ of $X$ and $Y$ equals $X\circ Y=Y\circ X$ if
$S(X,Y)=\emptyset$, and equals $\hat{1}$ otherwise.

An {\em affine oriented matroid} is a triple $(E,\L,g)$, where
$(E,\L)$ is an oriented matroid and $g\in E$ is a distinguished
element which is not a loop. Recall that $g$ is a {\em loop} if
$X_g=0$ for all $X\in \L$. We now define the bounded complex as in
Definition 4.5.1 of \cite{BLSWZ}. For an affine oriented matroid
$(E,\L,g)$ let
$$
\L^+ = \{X\in \L\ :\ X_g=+\}\ \ \hbox{and}\ \ \widehat{\L}^+=\L^+
\cup \{0,\hat{1}\}.
$$
With the induced order as a subset of $\widehat{\L}$, we call
$\widehat{\L}^+$ the {\em affine face lattice} of $(E,\L,g)$. The
{\em bounded complex} of $(E,\L,g)$ is
$$
\L^{++}=\{X\in \L^+\ :\ \L_{\leq X} \subseteq \widehat{\L}^+\}.
$$

Let $\T$ denote the set of maximal covectors (called {\em topes})
of $\L$. Let $B\in\T$. Then the {\em tope poset} $\T(\L,B)$ is a
partial order on the set $\T$ defined by $T\leq T'$ if and only if
$S(B,T) \subseteq S(B,T')$. The following lemma is easy to deduce
from this definition (see Corollary 4.2.11 of~\cite{BLSWZ}).

\begin{lemma}\label{same structure}
Let $T< T'$ in $\T(\L,B)$. Then the interval $[T,T']$ has the same
structure in $\T(\L,T)$ as in $\T(\L,B)$.
\end{lemma}

The following is a theorem of Lawrence (see Proposition 4.3.2
of~\cite{BLSWZ}).

\begin{proposition}\label{linear extension}
Every linear extension of $\T(\L,B)$ is a shelling. Therefore
$\L\setminus \{0\}$ is the face poset (with the empty face
excluded) of a shellable regular cell decomposition of an
$(r-1)$-sphere, where $r$ is the rank of the underlying matroid.
\end{proposition}

We also need the following more general result (see Corollary
4.3.7 of~\cite{BLSWZ}).

\begin{proposition}\label{interval}
Let $(X,Y)$ be an interval in $\widehat{\L}$. Then $(X,Y)$ is
isomorphic to the face poset of a shellable regular cell
decomposition of the $(\rank(Y)-\rank(X)-2)$-sphere.
\end{proposition}

$\L^+$ is an order filter in $\L$. It is pure of length $r-1$
where $r$ is the rank of the underlying matroid, i.e., every
maximal chain is of the same length $r-1$ (see Proposition 4.5.3
of~\cite{BLSWZ}). In particular, the minimal covectors in $\L^+$
are also atoms in $\L$. It is now easy to see that $\L^{++}$ is
never empty : if $X$ is a minimal covector in $\L^+$ then $X$ is
an atom of $\L$, therefore $X\in \L^{++}$.  In the realizable case
this corresponds to the fact that the bounded complex of an {\em
essential} hyperplane arrangement is nonempty. The poset $\L^{++}$
is an order ideal of $\L\setminus \{0\}$, hence the face poset of
a subcomplex of the $(r-1)$-sphere. For other basic facts about
$\L$, $\L^+$ and $\L^{++}$ we refer the reader to Chapter 4
of~\cite{BLSWZ}.

In this paper, the same symbol will often be used to denote a
regular cell complex and its face poset if no confusion can arise.
In particular we use $\L^{++}$ to denote both the poset and the
underlying cell complex. Some useful results about the bounded
complex from~\cite{D} are collected in the following proposition.

\begin{proposition}[\cite{D}]\label{pure}
Let $(E,\L,g)$ be an affine oriented matroid.
\begin{enumerate}
\item $\L^{++}$ is pure.

\item All the maximal covectors in the bounded complex $\L^{++}$
have the same support, say, $E_1$. Let $\L_1^{++}$ denote the
bounded complex of $(E_1, \L_1, g)$, where $\L_1$ is the deletion
$\L\backslash (E- E_1) = \{X|_{E_1} : X\in \L\}$. Then $\L_1^{++}
\cong \L^{++}$.
\end{enumerate}
\end{proposition}

\begin{proof}
Part (1) is Corollary 3.3 of~\cite{D}. Part (2) follows from
Theorem 3.2(2) and Theorem 5.1 of~\cite{D}.
\end{proof}

A matroid of rank $r$ is {\em uniform} if every $r$-element subset
of the ground set $E$ is a basis. An oriented matroid $(E,\L)$ is
{\em uniform} if its underlying matroid is uniform. Similarly an
affine oriented matroid $(E, \L, g)$ is {\em uniform} if $(E, \L)$
is uniform. The realizable uniform affine oriented matroids
correspond to exactly the simple hyperplane arrangements. The
following proposition characterizes uniform matroids in terms of
several different systems of axioms. We omit its straightforward
proof.

\begin{proposition}
Let $M$ be a matroid of rank $r$ on the ground set $E$. Then $M$
is uniform if and only if it satisfies one of the following
equivalent conditions:
\begin{enumerate}

\item Every $r$-element subset of $E$ is a basis.

\item $I\subseteq E$ is independent if and only if $|I|\leq r$.

\item $F\subseteq E$ is a flat if and only if either $F=E$ or
$|F|\leq r-1$.

\item For any subset $A\subseteq E$,
$$
\rank(A)=
\left\{%
\begin{array}{ll}
    |A|, & \hbox{if $|A|\leq r-1$;} \\
    r, & \hbox{otherwise.} \\
\end{array}%
\right.
$$
\end{enumerate}

\end{proposition}

For an oriented matroid $(E, \L)$, recall that the set $L=\{z(X) :
X\in \L\}$ is the collection of flats of the underlying matroid.
The map $z:\L\rightarrow L$ is a cover-preserving, order-reversing
surjection of $\L$ onto the geometric lattice $L$. Therefore we
have the following characterization of a uniform oriented matroid.

\begin{corollary}\label{uniform}
An oriented matroid $(E,\L)$ of rank $r$ is uniform if and only if
it satisfies one of the following equivalent conditions:
\begin{enumerate}
\item For every subset $F\subseteq E$ with $|F| \leq r-1$, there
exists a covector $X\in \L$ with $z(X)=F$.

\item $\rank(X) = r-|z(X)|$ for all $X \in \L \setminus \{0\}$.
\end{enumerate}
\end{corollary}

\section{Outline of the proof}

Let us review some necessary terminology from $PL$ topology.
Recall that a simplicial complex $K$ is a {\em $PL$ d-ball} if $K$
and the standard $d$-simplex have isomorphic subdivisions. A
simplicial complex $K$ is a {\em $PL$ $(d-1)$-sphere} if $K$ and
the boundary of the standard $d$-simplex have isomorphic
subdivisions. A simplicial complex $K$ is a {\em $PL$
$d$-manifold} if the link of every vertex is either a $PL$
$(d-1)$-sphere or a $PL$ $(d-1)$-ball.

Recall that the {\em order complex} of a poset $P$ is the abstract
simplicial complex whose vertices are the elements of $P$ and
whose faces are the chains $x_0<x_1<\cdots<x_k$ in $P$. The
geometric realization of the order complex will also be denoted by
$\Delta(P)$, or even just by $P$ if no confusion can arise.

Let $\Gamma$ be a {\em regular} cell complex with face poset
$\F(\Gamma)$. Then the order complex $\Delta(\F(\Gamma))$ is
homeomorphic to $\Gamma$ (see Proposition 4.7.8 of~\cite{BLSWZ}).
It is the first {\em barycentric subdivision} of $\Gamma$. The
$PL$ definitions are now applicable to regular cell complexes. A
regular cell complex $\Gamma$ is a $PL$ $d$-ball if and only if
its simplicial subdivision $\Delta(\F(\Gamma))$ is a $PL$
$d$-ball, and similarly for $PL$ spheres and $PL$ manifolds (see
Lemma 4.7.25 of~\cite{BLSWZ}).

The main objective of this paper is to prove the following
theorem.

\begin{theorem}\label{main}
Let $(E,\L, g)$ be a uniform affine oriented matroid. Then its
bounded complex $\L^{++}$ is a $PL$ ball.
\end{theorem}

The tool of {\em collapsing} will be used to prove the theorem.
Let $\Gamma$ be a regular cell complex, and suppose that
$\sigma\in\Gamma$ is a proper face of exactly one face
$\tau\in\Gamma$. Then the complex
$\Gamma'=\Gamma\setminus\{\sigma,\tau\}$ is obtained from $\Gamma$
by an {\em elementary collapse}. Note that the condition on
$\sigma$ and $\tau$ implies that $\tau$ is a maximal face of
$\Gamma$ and $\sigma$ is a maximal proper face of $\tau$. If
$\Gamma$ can be reduced to a single point by a sequence of
elementary collapses, then $\Gamma$ is {\em collapsible}. The
following is Theorem 6.11 of~\cite{D}.

\begin{theorem}[\cite{D}]\label{collapsible}
Let $(E,\L, g)$ be an affine oriented matroid. Then its bounded
complex $\L^{++}$ is collapsible.
\end{theorem}

\begin{remark}\rm
The collapsing in $PL$ topology is more restrictive than our above
definition. To distinguish between the two let us define {\em $PL$
collapsing} here. We follow the notation in~\cite{RS}. Suppose
that $X\supset Y$ are {\em locally conical sets} (called {\em
polyhedra} in~\cite{RS}), $B^n$ is a $PL$ $n$-ball, and $B^{n-1}$
is a $PL$ $(n-1)$-ball contained in the boundary of $B^n$. If
$X=Y\cup B^n$ and $Y\cap B^n = B^{n-1}$, then we say that there is
an {\em elementary $PL$ collapse} of $X$ on $Y$. We say that $X$
is $PL$ collapsible if $X$ reduces to a point via a sequence of
elementary $PL$ collapses.
\end{remark}

We want to show that $\L^{++}$ is $PL$ collapsible. For this we
need the following proposition (see Proposition 4.7.26
of~\cite{BLSWZ}).

\begin{proposition}\label{sphere}
Let $\Gamma$ be a regular cell decomposition of the $d$-sphere. If
$\Gamma$ is shellable then $\Gamma$ is a $PL$ sphere. If $\Gamma$
is a $PL$ sphere then every closed cell in $\Gamma$ is a $PL$
ball.
\end{proposition}

\begin{corollary}
$\L^{++}$ is $PL$ collapsible.
\end{corollary}

\begin{proof}
Regular cell complexes are locally conical sets since they admit
simplicial subdivisions. The cells in $\L^{++}$ are all $PL$ balls
by Proposition~\ref{linear extension} and
Proposition~\ref{sphere}. Therefore the bounded complex $\L^{++}$
is in fact $PL$ collapsible.
\end{proof}

We need the following fact from $PL$ topology. For a proof see
Corollary 3.28 of~\cite{RS}.

\begin{theorem}
If a $PL$ manifold is $PL$ collapsible, then it is a $PL$ ball.
\end{theorem}

To prove the main Theorem~\ref{main}, it remains to show that the
following is true.

\begin{lemma}\label{key}
Let $(E,\L, g)$ be a uniform affine oriented matroid. Then its
bounded complex $\L^{++}$ is a $PL$ manifold.
\end{lemma}

This lemma will be proven in the next section, by showing that the
link of a vertex in $\Delta(\L^{++})$ is either a $PL$ ball or a
$PL$ sphere.

\section{The links in $\Delta(\L^{++})$}

Let $(E,\L,g)$ be an affine oriented matroid (not necessarily
uniform). Let $X\in \L^{++}$ be a covector in its bounded complex.
For simplicity we will not distinguish $\{X\}$ and $X$. Then $X$
is a vertex in the order complex $\Delta(\L^{++})$. The link of
$X$ in $\Delta(\L^{++})$ will be denoted by
$$
\link(X,\Delta(\L^{++})) := \{\sigma\in \Delta(\L^{++}):\, X\notin
\sigma\ \hbox{and}\ \{X\}\cup\sigma \in \Delta(\L^{++})\}.
$$

Recall that the {\em join} of two simplicial complexes $K_1$ and
$K_2$ on {\em disjoint} vertex sets is
$$
K_1*K_2 := \{ \sigma_1\cup \sigma_2:\, \sigma_1\in K_1,
\sigma_2\in K_2\}.
$$
Let $P$ be a finite poset and $x\in P$. Then it is easy to see
that
$$
\link(x,\Delta(P))= \Delta(P_{<x}) * \Delta(P_{>x}),
$$
where $P_{<x}=\{y\in P: y<x\}$ and similarly for $P_{>x}$. It
follows that
$$
\link(X,\Delta(\L^{++}))= \Delta(\L^{++}_{<X}) *
\Delta(\L^{++}_{>X}).
$$

Let us first consider $\Delta(\L^{++}_{<X})$. Recall that
$\L^{++}$ is an order ideal in $\L\setminus \{0\}$. Hence
$\L^{++}_{<X}$ is the same as the interval $(0, X)$ in $\L$. By
Proposition~\ref{interval}, $(0,X)$ is isomorphic to the face
poset of a shellable regular cell decomposition of a sphere of
dimension $\rank(X)-2$. Therefore $\L^{++}_{<X}$ is a $PL$ sphere
by Proposition~\ref{sphere}, and so is $\Delta(\L^{++}_{<X})$.
Note that when $X$ is of rank one, $\Delta(\L^{++}_{<X})$ is the
complex $\{\emptyset\}$ which we consider as a sphere of dimension
$-1$.

Next let us consider $\Delta(\L^{++}_{>X})$. First we make some
reductions. By Proposition~\ref{pure}(2) we may assume that
$\L^{++}$ is full dimensional, so that the maximal covectors in
$\L^{++}$ are topes of $\L$. If $\L^{++}_{>X}=\emptyset$ then
$\Delta(\L^{++}_{>X})$ is a sphere of dimension $-1$. If
$\L^{++}_{>X}=\L_{>X}$ then $\Delta(\L^{++}_{>X})$ is the first
barycentric subdivision of $\L_{>X}$ which is a shellable regular
cell decomposition of a sphere, hence $\Delta(\L^{++}_{>X})$ is a
$PL$ sphere. In what follows we assume that  $\emptyset \subsetneq
\L^{++}_{>X} \subsetneq \L_{>X}$.

From now on we shall use the fact that $(E,\L,g)$ is uniform. The
following several lemmas are not true for general affine oriented
matroids. They are the crux of our proof of the conjecture.

\begin{lemma}\label{cross}
Let $(E,\L)$ be a uniform oriented matroid and $X\in \L\setminus
\{0\}$. Then the map
\begin{align*}
    d:\L_{\geq X} &\rightarrow \{+,-,0\}^{z(X)}\\
    Y &\mapsto Y\backslash \supp(X)
\end{align*}
is an isomorphism of posets.
\end{lemma}

\begin{proof}
Let $r$ be the rank of $(E,\L)$. By Corollary~\ref{uniform},
$|z(X)|\leq r-1$ since $X\neq 0$. Again by
Corollary~\ref{uniform}, for every $e\in z(X)$ there exists a
covector $Y\in\L$ with $z(Y) = z(X) -\{e\}$. It follows that for
every $e\in z(X)$ there is a covector $Z\in \L\backslash \supp(X)$
with $\supp(Z)=\{e\}$. The covector axioms then imply that
$\L\backslash \supp(X)= \{+,-,0\}^{z(X)}$. Finally note that
$\L_{\geq X}$ is isomorphic to the deletion $\L\backslash
\supp(X)$ via the map $Y \mapsto Y\backslash \supp(X)$.
\end{proof}

Since the same symbol is often used to denote a face poset and its
underlying regular cell complex in this paper, when a face poset
$P$ is said to be simplicial or shellable it is meant that the
underlying complex is simplicial or shellable.

\begin{corollary}\label{simplicial}
Let $(E,\L)$ be a uniform oriented matroid and $X\in \L\setminus
\{0\}$. Then $\L_{>X}$ is simplicial.
\end{corollary}

\begin{proof}
By Lemma~\ref{cross}, $\L_{>X}$ is isomorphic to the face poset of
the boundary of a cross polytope, which is simplicial.
\end{proof}

Recall that without loss of generality $\L^{++}$ is assumed to be
full dimensional. We may also assume that $|E|>1$ to exclude the
trivial case $E=\{g\}$. Under these assumptions, for every $X\in
\L^{+}$, the deletion $X\backslash g \neq 0$. Otherwise, if
$X\backslash g = 0$, then it is easy to show that $\L^{++}
=\{X\}$. This contradicts the full-dimensionality of $\L^{++}$.

\begin{corollary}\label{boundary}
If $X\in \L^+$, then $X\backslash g \in \L/g$ if and only if
$X\notin \L^{++}$.
\end{corollary}

\begin{proof}
Let $Z\in\{+,-,0\}^E$ be the sign vector defined by $Z_g = 0$ and
$Z_e = X_e$ otherwise. Note that $Z \neq 0$ since $X\backslash g
\neq 0$.

If $X\notin \L^{++}$, then there exists $0<Y<X$ such that $Y_g
=0$. Note that $Y\leq Z$ as sign vectors. Applying
Lemma~\ref{cross} to $\L_{\geq Y}$, we see that $Z$ is a covector.
Therefore $X\backslash g \in \L/g$. Conversely, if $X\backslash g
\in \L/g$, then $Z\in \L$. It follows that $X\notin \L^{++}$ since
$0<Z<X$ and $Z_g = 0$.
\end{proof}

For $X\in \L^{++}$ let $\C_X$ denote the set of topes that are in
$\L_{>X}=\L^+_{>X}$ but not in $\L^{++}_{>X}$. Let $\D_X$ denote
the set of topes $T$ of the contraction $\L/g$ with the following
property: $T\in \D_X$ if and only if $T_e=X_e$ for all $e\in
\supp(X\backslash g)$. Equivalently, $T\in \D_X$ if and only if
$T\geq X\backslash g$. We should, however, note that $X\backslash
g \notin \L/g$ in this case.

\begin{lemma}\label{bijection}
A tope $T\in \L_{>X}$ is in $\C_X$ if and only if the deletion
$T\backslash g$ is a tope in $\D_X$. Moreover the map $r:
\C_X\rightarrow \D_X$ defined by $r(T) = T\backslash g$ is a
bijection.
\end{lemma}

\begin{proof}

The first statement follows from Corollary~\ref{boundary}. To
prove the second statement, let $T\in\D_X$. Define $h(T)
\in\{+,-,0\}^E$ by $h(T)_g=X_g=+$ and $h(T)_e=T_e$ for all $e\in
E\setminus \{g\}$. Define $i(T)\in\{+,-,0\}^E$ by $i(T)_g=0$ and
$i(T)_e=T_e$ for all $e\in E\setminus \{g\}$. Then $i(T) \in \L$
since $T\in \L/g$. Therefore $h(T)= i(T) \circ X \in \L^{+}$.
Moreover $h(T) \notin \L^{++}$ since $0<i(T)<h(T)$ in $\L$. Hence
$h(T)\in\C_X$. Note that $h$ is the inverse map of $r$, showing
that $r$ is bijective.
\end{proof}

Let $[\D_X] = \{Y\in \L/g: 0<Y \leq T\ \hbox{for some}\ T\in \D_X
\}$. Then $[\D_X]$ is a subcomplex of the shellable sphere $\L/g$.

\begin{lemma}
$[\D_X]$ is shellable.
\end{lemma}
\begin{proof}
Fix a tope $B\in \D_X$. For a tope $T'\in\D_X$ we have
$T'_e=B_e=X_e$ for all $e\in\supp(X\backslash g)$. If $T\in\L/g$
and $T\leq T'$ in the tope poset $\T(\L,B)$, then $S(B,T)\subseteq
S(B,T')$. Hence $T_e=T'_e=B_e=X_e$ for all $e\in\supp(X\backslash
g)$, so $T\in\D_X$. Therefore $\D_X$ is an order ideal in the tope
poset, so there is a linear extension of the tope poset in which
$\D_X$ is an initial segment. Hence $\D_X$ is an initial segment
of a shelling of the topes in $\L/g$ by Proposition~\ref{linear
extension}, proving that $[\D_X]$ is shellable.
\end{proof}

Let $[\C_X] = \{Y \in \L_{>X}: Y \leq T\ \hbox{for some}\ T\in
\C_X \}$. It turns out that a shelling on $\D_X$ induces a
shelling on $\C_X$. The shellability of a regular cell complex is
defined recursively (see Definitions 4.7.14 and 4.7.17
of~\cite{BLSWZ}). However in the case of a simplicial complex, the
definition of a shelling is much simpler. Fortunately for us, the
poset $[\C_X]$ is an order ideal in $\L_{>X}$ and hence the face
poset of a simplicial complex by Lemma~\ref{simplicial}.

First let us formulate the definition of a shellable simplicial
complex in terms of its face poset. Let $P$ be the face poset of a
pure $d$-dimensional simplicial complex. Let $\widehat{P} = P\cup
\{\hat{0},\hat{1}\}$ denote the (augmented) face lattice by adding
new elements such that $\hat{0}<p<\hat{1}$ for all $p\in P$. A
linear ordering $c_1, c_2, \ldots, c_t$ of the coatoms of
$\widehat{P}$ is a {\em shelling} if for all $1\leq i\leq j\leq t$
there exists $1\leq k< j$ such that $c_i\wedge c_j \leq c_k\wedge
c_j \lessdot c_j$, where $\lessdot$ is the covering relation in
the poset. Note that if $P$ is the face poset of a regular cell
complex whose augmented face poset is a lattice (as in the case of
$[\D_X]$), then the above condition on the linear order of coatoms
is necessary, although not sufficient, for being a shelling.

\begin{lemma}\label{shellable}
$[\C_X]$ is shellable.
\end{lemma}

\begin{proof}
Fix a tope $B \in \D_X$. Let $d_1 (=B), d_2, \ldots, d_t$ be the
shelling of $\D_X$ obtained from an initial segment of a linear
extension of the tope poset $\T(\L/g, B)$. Let $c_i = h(d_i) \in
\C$ for all $1\leq i\leq t$, where $h$ is as defined in the proof
of Lemma~\ref{bijection}. We want to show that $c_1, c_2, \ldots,
c_t$ is a shelling of $[\C_X]$. As in the above definition of
shelling we shall work with the augmented face lattices
$\widehat{[\C_X]}$ and $\widehat{[\D_X]}$. We identify $\hat{0}$
in $\widehat{[\C_X]}$ with $X$, and $\hat{0}$ in
$\widehat{[\D_X]}$ with $0$ in $\L/g$. Let $1\leq i< j\leq t$, we
want to find $1\leq k <j$ such that $c_i\wedge c_j \leq c_k\wedge
c_j \lessdot c_j$.

{\em Case 1}: If $c_i\wedge c_j \notin \L^{++}$, then $(c_i\wedge
c_j)\backslash g \in \L/g$ by Corollary~\ref{boundary}. Therefore
$(c_i\wedge c_j)\backslash g \leq d_i \wedge d_j$. By the shelling
of $\D_X$ there exists $1\leq k <j$ such that $d_i\wedge d_j \leq
d_k \wedge d_j \lessdot d_j$. Since $(c_i\wedge c_j)\backslash g
\leq d_i\wedge d_j \leq d_k$, we get $c_i\wedge c_j\leq c_k$.
Hence $c_i \wedge c_j \leq c_k \wedge c_j$. Since $\L/g$ is
uniform we have $|z(d_k\wedge d_j)|=1$. It follows that
$|z(c_k\wedge c_j)|=1$ and hence $c_k\wedge c_j \lessdot c_j$.

{\em Case 2}: If $c_i \wedge c_j \in \L^{++}$, then consider
$\D_{c_i\wedge c_j}$ which is a subset of $\D_X$. The interval
$[d_i, d_j]$ has the same structure in $\T(\L/g, d_i)$ and
$\T(\L/g, B)$ by Lemma~\ref{same structure}. Note that $[d_i, d_j]
\subseteq \D_{c_i\wedge c_j}\subseteq \D_X$ since $\D_{c_i\wedge
c_j}$ is an order ideal in the tope poset $\T(\L/g, d_i)$. There
is a linear extension of $\T(\L/g, d_i)$ (hence a shelling of
$\L/g$) such that
\begin{enumerate}
\item $[d_i, d_j]$ is an initial segment;

\item the linear order of the elements in $[d_i, d_j]$ is the same
as the restriction of the shelling of $\D_X$ on $[d_i, d_j]$.
\end{enumerate}
Therefore there exists $d_k \in [d_i, d_j]$ (so $i\leq k < j$)
such that $d_k\wedge d_j \lessdot d_j$. Once again, since $\L/g$
is uniform we have $|z(d_k\wedge d_j)| = 1$. It follows that
$|z(c_k\wedge c_j)|=1$ and hence $c_k\wedge c_j \lessdot c_j$.
Finally $d_k \in \D_{c_i\wedge c_j}$ implies that $c_i \wedge c_j
\leq c_k$, so $c_i\wedge c_j \leq c_k\wedge c_j$.

\end{proof}

\begin{corollary}
$[\C_X]$ is a $PL$ ball.
\end{corollary}

\begin{proof}
Since $[\C_X]$ is shellable, it is either a $PL$ ball or a $PL$
sphere. Since $\L^{++}_{>X}$ is nonempty, $[\C_X]$ is a {\em
proper} subset of $\L_{>X}$ which is a $PL$ sphere of the same
dimension. Therefore $[\C_X]$ has to be a $PL$ ball.
\end{proof}

The following lemma is known as Newman's Theorem (see Theorem
4.7.21(iii) of~\cite{BLSWZ}).

\begin{lemma}
The closure of the complement of a $PL$ $d$-ball embedded in a
$PL$ $d$-sphere is itself a $PL$ $d$-ball.
\end{lemma}

$[\C_X]$ is a full dimensional $PL$ ball embedded in the $PL$
sphere $\L_{>X}$. The closure of its complement is exactly
$\L^{++}_{>X}$. It follows that $\L^{++}_{>X}$ is a $PL$ ball, and
so is $\Delta(\L^{++}_{>X})$.

We have seen that $\Delta(\L^{++}_{<X})$ is a $PL$ sphere and
$\Delta(\L^{++}_{>X})$ is either a $PL$ sphere or a $PL$ ball, so
$$
\link(X,\Delta(\L^{++}))= \Delta(\L^{++}_{<X})
* \Delta(\L^{++}_{>X})
$$
is either a $PL$ sphere or a $PL$ ball by the following lemma (see
Proposition 2.23 of~\cite{RS}).

\begin{lemma}
The join of two $PL$ spheres is a $PL$ sphere. The join of a $PL$
sphere and a $PL$ ball is a $PL$ ball.
\end{lemma}

By Proposition~\ref{pure}(1) $\L^{++}$ is pure of dimension $r-1$,
so $\link(X,\Delta(\L^{++}))$ is either a $PL$ $(r-2)$-sphere or a
$PL$ $(r-2)$-ball. We conclude that $\Delta(\L^{++})$ is a $PL$
$(r-1)$-manifold, and so is $\L^{++}$. Since $\L^{++}$ is known to
be $PL$ collapsible, it is in fact a $PL$ $(r-1)$-ball. The proof
of our main Theorem~\ref{main}  is now complete.

\section{Final remarks}

Zaslavsky~\cite{Za} in fact also conjectures that the bounded
complex of a hyperplane arrangement is a ball as long as there is
no parallelism among the hyperplanes and intersections.
Generalizing to affine oriented matroids, we get the following
conjecture.

\begin{conjecture}
Let $(E,\L,g)$ be an affine oriented matroid. If the contraction
$\L/g$ is uniform, then the bounded complex $\L^{++}$ is a ball.
\end{conjecture}

Another open question in~\cite{BLSWZ} asks whether the bounded
complex of a simplicial affine oriented matroid is a ball. We
phrase it as a conjecture.

\begin{conjecture}
Let $(E,\L,g)$ be a simplicial affine oriented matroid. Then the
bounded complex $\L^{++}$ is a ball.
\end{conjecture}

More generally one can ask for which kind of affine oriented
matroids the bounded complex is a ball. A related question is for
which kind of affine oriented matroids the bounded complex is
shellable. If the bounded complex is shellable then it must be a
ball, but not vice versa. It is not known whether the bounded
complex of a uniform affine oriented matroid is shellable.

\newcommand{\journalname}[1]{\textrm{#1}}
\newcommand{\booktitle}[1]{\textrm{#1}}


\begin{thebibliography}{9}

\bibitem{BZ}
A. Bj\"orner and G. M. Ziegler, \textit{Shellability of oriented
matroids}, Abstract, Workshop on Simplicial Complexes, IMA, Univ.
of Minnesota, March 1988; and Abstract, Conf. on Ordered Sets,
Oberwolfach, April 1988.

\bibitem{BLSWZ}
A. Bj\"orner, M. Las Vergnas, B. Sturmfels, N. White and G.
Ziegler, \booktitle{Oriented Matroids}, 2nd Ed., Cambridge Univ.
Press, 1999.

\bibitem{D}
X. Dong, \textit{On the bounded complex of an affine oriented
matroid}, \journalname{Discrete Comput. Geom.}, {\bf 35}(2006),
no. 3, 457--471.

\bibitem{L}
M. Las Vergnas, \textit{Matr\"oides orientables}, \journalname{C.
R. Acad. Sci. Paris S\'er. A-B}, {\bf 280} (1975), Ai, A61--A64.

\bibitem{Ma}
B. Mazur, \textit{A note on some contractible 4-manifolds},
\journalname{Annals of math.}, {\bf 73} (1961), no. 1, 221--228.

\bibitem{RS}
C. Rourke and B. Sanderson, \booktitle{Introduction to
piecewise-linear topology}, Springer-Verlag, 1982.

\bibitem{St}
R. Stanley, \booktitle{An introduction to hyperplane
arrangements}, Lecture notes, IAS/Park City Mathematics Institute,
2004.

\bibitem{Za}
T. Zaslavsky, \textit{Facing up to arrangements: face-count
formulas for partitions of space by hyperplanes},
\journalname{Mem. Amer. Math. Soc.} {\bf 1} (1975), issue 1, no.
154.

\bibitem{Zi}
G. M. Ziegler, \textit{The face lattice of hyperplane
arrangements}, \journalname{Discrete Math.}, {\bf 74} (1988),
233--238.








\end{thebibliography}
\end{document}